\documentclass[a4paper]{amsart}

\usepackage{amscd}
\usepackage[T1]{fontenc}
\usepackage[latin1]{inputenc}
\usepackage[pdfauthor={Stefan Kebekus, Luis Solá and Matei Toma},
pdftitle={Rationally connected foliations after Bogomolov and McQuillan},
pdfkeywords={algebraic variety, foliation, rationally connected variety},
pdfview={Fit}]{hyperref}
\usepackage[arrow,curve,matrix]{xy}
\usepackage{times}
\usepackage{epsfig}
\newcommand{\sF}{\mathcal{F}}

\renewcommand{\P}{\mathbb{P}}

\renewcommand{\O}{\mathcal{O}}
\newcommand{\Q}{\mathbb{Q}}

\DeclareMathOperator{\charac}{char}
\DeclareMathOperator{\Chow}{Chow}
\DeclareMathOperator{\codim}{codim}
\DeclareMathOperator{\Hilb}{Hilb}
\DeclareMathOperator{\Hom}{Hom}

\DeclareMathOperator{\length}{length}

\DeclareMathOperator{\Pic}{Pic}
\DeclareMathOperator{\red}{red}
\DeclareMathOperator{\reg}{reg}

\DeclareMathOperator{\rank}{rank}

\DeclareMathOperator{\Spec}{Spec}

\DeclareMathOperator{\id}{id}

\newcommand{\ilabel}[1]{\newcounter{#1}\setcounter{#1}{\value{enumi}}}
\newcommand{\iref}[1]{\setcounter{enumi}{\value{#1}}\labelenumi}

 \theoremstyle{plain}
 \newtheorem{thm}{Theorem}
 \newtheorem{defn}[thm]{Definition}

 \numberwithin{figure}{section}

 \theoremstyle{plain}
 \newtheorem{cor}[thm]{Corollary}
 \newtheorem{lem}[thm]{Lemma}

 \theoremstyle{plain}
 \newtheorem{prop}[thm]{Proposition}

 \theoremstyle{remark}
 \newtheorem{fact}[thm]{Fact}
 \newtheorem{rem}[thm]{Remark}
 
\theoremstyle{remark}
 \newtheorem{notation}[thm]{Notation}
 \newtheorem{assumption}[thm]{Assumption}
 
 \newtheorem{question}[thm]{Question}

\def\factor#1.#2.{\left. \raise 2pt\hbox{$#1$} \right/\hskip -2pt\raise -2pt\hbox{$#2$}}

\begin{document}

\title{Rationally connected foliations after Bogomolov and McQuillan}

\date{\today} 

\author[S.~Kebekus, L.~Solá Conde, and M.~Toma] {Stefan~Kebekus, Luis~Solá Conde and Matei Toma}

\address{Stefan Kebekus, Mathematisches Institut, Universität  zu Köln,
  Weyertal 86-90, 50931 Köln, Germany}
\email{stefan.kebekus@math.uni-koeln.de}
\urladdr{http://www.mi.uni-koeln.de/\~{}kebekus}

\address{Luis Solá Conde, Mathematisches Institut, Universität zu Köln,
  Weyertal 86-90, 50931 Köln, Germany} \email{lsola@math.uni-koeln.de}

\address{Matei Toma, Mathematisches Institut, Universität zu Köln,
  Weyertal 86-90, 50931 Köln, Germany and Mathematical Institute of the
  Romanian Academy, Bucharest} \email{matei@math.uni-koeln.de}

\thanks{The authors thank their hosting institution, Universit\"at zu K\"oln.
  The first two authors were supported in full or in part by the
  Forschungsschwerpunkt ``Globale Methoden in der komplexen Analysis''
  of the Deutsche Forschungsgemeinschaft. A part of this paper was
  worked out while Stefan Kebekus visited the Korea Institute for
  Advanced Study.  He would like to thank Jun-Muk Hwang for the
  invitation.}

\begin{abstract}
  This paper is concerned with sufficient criteria to guarantee that a
  given foliation on a normal variety has algebraic and rationally
  connected leaves. Following ideas from a preprint of
  Bogomolov-McQuillan and using the recent work of
  Graber-Harris-Starr, we give a clean, short and simple proof of
  previous results. Apart from a new vanishing theorem for vector
  bundles in positive characteristic, our proof employs only standard
  techniques of Mori theory and does not make any reference to the
  more involved properties of foliations in characteristic $p$.
  
  We also give a new sufficient condition to ensure that all leaves
  are algebraic.
  
  The results are then applied to show that $\mathbb Q$-Fano varieties
  with unstable tangent bundles always admit a sequence of partial
  rational quotients naturally associated to the Harder-Narasimhan
  filtration.
\end{abstract}

\maketitle

\setcounter{tocdepth}{1}
\tableofcontents

\section{Introduction}

Since Ekedahl's work on foliations in positive characteristic and
Miyaoka's landmark paper \cite{Miy85}, foliations of algebraic
varieties have met with considerable interest from algebraic geometers
and number theorists alike.

This paper is concerned with a sufficient criterion to guarantee that
a given singular foliation on a normal variety has algebraic and
rationally connected leaves. More precisely, using a vanishing theorem
for vector bundles on curves in characteristic $p$, we give a simple
proof of the following result.

\begin{thm}[cf.~\protect{\cite[thm.~0.1]{BMcQ}}]\label{thm:algleafs}
  Let $X$ be a normal complex projective variety, $C \subset X$ a
  complete curve which is entirely contained in the smooth locus
  $X_{\reg}$, and $\sF \subset T_X$ a (possibly singular) foliation
  which is regular along $C$. Assume that the restriction $\sF|_C$ is
  an ample vector bundle on $C$. If $x \in C$ is any point, the leaf
  through $x$ is algebraic. If $x \in C$ is general, the closure of
  the leaf is rationally connected.
\end{thm}

Bogomolov and McQuillan state a stronger result where the foliation is
allowed to have singularities at points of $C$. Since this paper
strives for simplicity, not completeness, we do not consider this
extra complication here.

On the other hand, if $\sF$ is non-singular, we can use the Reeb
stability theorem to prove the algebraicity of all leaves.

\begin{thm}\label{thm:reeb}
  In the setup of Theorem~\ref{thm:algleafs}, if $\sF$ is regular,
  then all leaves are rationally connected submanifolds.
\end{thm}

In fact, a stronger statement holds, see Theorem~\ref{thm:reeb2} in
section~\ref{sec:reeb}.

As an immediate corollary to Theorem~\ref{thm:algleafs}, we prove a
refinement of Miyaoka's characterization of uniruled varieties. Before
stating this corollary, we need to introduce the following notation.

\begin{notation}
  Let $X$ be a normal projective variety, and let $q : X \dasharrow Q$
  be the rationally connected quotient, defined through the maximally
  rationally connected fibration (``MRC-fibration'') of a
  desingularization of $X$, cf.~\cite[IV. 5.3, 5.5]{K96}. Further,
  suppose that $C \subset X$ is a subvariety which is not contained in
  the singular locus of $X$, and not contained in the indeterminacy
  locus of $q$, and that $\sF \subset T_X|_C$ is a subsheaf of the
  restriction of the tangent sheaf to $C$. We say that $\sF$ is
  \emph{vertical with respect to the rationally connected quotient},
  if $\sF$ is contained in $T_{X|Q}$ at the general point of $C$.
\end{notation}

\begin{notation}
  If $X$ is normal, we consider \emph{general complete intersection
    curves} in the sense of Mehta-Ramanathan, $C \subset X$. These are
  reduced, irreducible curves of the form $C = H_1 \cap \cdots \cap
  H_{\dim X-1}$, where the $H_i \in |m_i \cdot L_i|$ are general, the
  $L_i \in \Pic(X)$ are ample and the $m_i \in \mathbb N$ large
  enough, so that the Harder-Narasimhan filtration of $T_X$ commutes
  with restriction to $C$.  If the $L_i$ are chosen, we also call $C$
  a \emph{general complete intersection curve with respect to $(L_1,
    \ldots, L_{\dim X-1})$.}
  
  We refer to \cite{Flenner84} and \cite{Langer04} for a discussion
  and an explicit bound for the $m_i$.
\end{notation}

\begin{cor}\label{cor:ratconnleaf}
  Let $X$ be a normal complex-projective variety and $C \subset X$ a
  general complete intersection curve. Assume that the restriction
  $T_X|_C$ contains an ample locally free subsheaf $\sF_C$. Then
  $\sF_C$ is vertical with respect to the rationally connected
  quotient of $X$.
\end{cor}

This statement appeared first implicitly in
\cite[chap.~9]{SecondAsterisque}, but see
Remark~\ref{rem:problemsinsecondasterisque}. To our best knowledge,
the argument presented here gives the first complete proof of this
important result.

\medskip

We would like to thank Thomas Peternell for discussions and the
anonymous referee for careful reading and a number of helpful
comments. In particular, we owe the referee an elementary proof of our
vanishing result, Proposition~\ref{rem:atft_gets_pos}; our previous
approach relied on more involved results of Langer, \cite{Langer04}.
The statement of Theorem~\ref{thm:reeb} was suggested to us by János
Kollár. We would like to thank him very much for the discussion.

\subsection{Outline of the paper}

All our results here are principally based on a vanishing theorem for
vector bundles in finite characteristic,
Proposition~\ref{prop:vanishing} in Section~\ref{sec:VBAC} below. In
Section~\ref{sec:deformation} we have gathered a number of standard
facts about the space of relative deformations of morphisms.

With these preparations at hand, the proof of
Theorem~\ref{thm:algleafs}, which we give in
Section~\ref{sec:proofOf1} becomes reasonably short and quite
intuitive. The line of argumentation follows \cite{BMcQ}, but differs
in one important aspect: we employ Mori's standard method of
``reduction modulo characteristic $p$'', but using our vanishing
result, we do not have to use the partial Frobenius morphism for
foliated varieties that was introduced by Miyaoka. We believe that
this, and our use of the Graber-Harris-Starr result makes the proof
much cleaner and easier to understand.

Theorem~\ref{thm:reeb} is proven in section~\ref{sec:reeb} using Reeb
stability and standard facts on rationally connected submanifolds.
After this paper was accepted for publication we learned that in his
preprint \cite{Hoe05}, Andreas Höring has independently obtained
similar results in the case of regular foliations.

Corollary~\ref{cor:ratconnleaf}, which we prove in
Section~\ref{sec:cor2}, follows from Theorem~\ref{thm:algleafs} and
the fact that vector bundles on curves always contain a maximally
ample subbundle which appears in the Harder-Narasimhan filtration.
While the latter is possibly known among experts, we could not find
any reference whatsoever, and chose to include a full and detailed
proof.

The methods of Section~\ref{sec:cor2} allow us to find a sequence of
partial rationally connected quotients naturally associated to the
Harder-Narasimhan filtration of the tangent bundle. For a precise
statement, see Section~\ref{sec:interpret}.

\medskip

In the course of proving Theorem~\ref{thm:algleafs} and
Corollary~\ref{cor:ratconnleaf}, we liberally use ideas and methods
introduced elsewhere. We have therefore, in Sections~\ref{sec:att1}
and \ref{sec:att2}, gathered all the references we are aware of.

 \section{A vanishing criterion for vector bundles in positive characteristic}
 \label{sec:VBAC}
 
 Throughout the present section, let $C$ be a smooth curve of genus
 $g$ defined over an algebraically closed field $k$ of characteristic
 $p>0$. Further, let $E$ be a vector bundle on $C$. Even if $E$ is
 ample, it is of course generally false that $H^1(C, E) = 0$ or that
 $E$ has any sections. However, in Proposition~\ref{prop:vanishing} we
 will give a criterion to guarantee that the pull-back of $E$ via the
 Frobenius morphism is globally generated and satisfies a strong
 vanishing statement, even if the pull-back is twisted with certain
 ideal sheaves.
 
 We will use $\Q$-twists of vector bundles as presented in
 \cite[II,~6.2]{Laz04}, to which we refer for details. In our case, we
 identify rational numbers $\delta$ with numerical classes
 $\delta\cdot[P]\in N^1_{\Q}(C)$, where $P$ is a point in $C$. For
 every $\delta\in\Q$, the {\it $\Q$-twist} $E\langle\delta\rangle$ is
 defined as the ordered pair of $E$ and $\delta$. A $\Q$-twisted
 vector bundle is said to be {\it ample} if the class $c_1
 \bigl(\O_{\mathbb P_C(E)}(1) \bigr)+\pi^*(\delta)$ is ample on the
 projectivized bundle $\P_C(E)$, where $\pi$ denotes the natural
 projection. One defines the degree
 $\deg(E\langle\delta\rangle):=\deg(E)+\rank(E)\delta$. A quotient of
 $E\langle\delta\rangle$ is a $\Q$-twisted vector bundle of the form
 $E'\langle\delta\rangle$ where $E'$ is a quotient of $E$.  Pull-backs
 of $\Q$-twisted vector bundles are defined in the obvious way.
 
 In order to state the main result of this section, we introduce the
 following notation.

 \begin{notation}
   Let $F:C[1]\rightarrow C$ be the $k$-linear Frobenius morphism.
 \end{notation}
 
 \begin{defn}[Vanishing threshold]\label{def:atopft}
   Given a rational number $\delta$ we define
   $$
   b_p(\delta):=p\delta-2g+1, 
   $$
   where $g$ denotes the genus of $C$.
 \end{defn}

 The following trivial observation will later become important in the
 applications, when $C$ and $E$ are reductions modulo $p$ of objects
 that were initially defined in characteristic 0.

 \begin{rem}\label{rem:atft_gets_pos}
   If we view the number $b_p(\delta)$ purely numerically as a function of
   integers $g$, $p$ and $\delta$, then
   $$
   \lim_{p\to \infty} \frac{p}{b_p(\delta)} = \lim_{p\to \infty}
   \frac{p}{\lfloor b_p(\delta) \rfloor} = \dfrac{1}{\delta},
   $$
   where $\lfloor . \rfloor$ means ``round-down''.
 \end{rem}

 \begin{prop}[Vanishing after pull-back via Frobenius]\label{prop:vanishing}
   Let $E$ be a vector bundle of rank $r$ over $C$, and $\delta$ a
   positive rational number. Assume that $E\langle -\delta\rangle$ is
   ample and that the vanishing threshold $b_p(\delta)$ is
   non-negative. Then for every subscheme $B\subset C[1]$ of length
   smaller than or equal to $b_p(\delta)$ we have
   \begin{equation}
     \label{eq:van}
     H^1 \bigl(C[1], F^*(E)\otimes I_B \bigr) = \{ 0\}.
   \end{equation}
   Further, $F^*(E)\otimes I_B$ is globally generated.
 \end{prop}

 \begin{proof}
   By Serre duality, the vanishing statement~\eqref{eq:van} is
   equivalent to
   \begin{equation}
     \label{eq:nonexmor}
     \Hom_{C} \bigl( F^*(E),\, \O_{C[1]}(B) \otimes \omega_{C[1]}\bigr) = \{0\}.    
   \end{equation}
   Since $E\langle -\delta\rangle$ is ample, its pull-back by $F$ is
   ample, too. Therefore any quotient of $F^*(E\langle -\delta\rangle)
   = F^*(E)\langle -p\delta \rangle$ is ample.  In particular, any
   quotient of rank $1$ of $F^*(E)$ has degree larger than $p\delta$.
   It follows that (\ref{eq:nonexmor}) is fulfilled whenever
   $$
   \deg(\O_{C[1]}(B)\otimes\omega_{C[1]})\leq p\delta,
   $$
   that is $\#(B)\leq b_p(\delta)+1$.

   To show that $F^*(E)\otimes I_B$ is globally generated, observe
   that, using the above argument, for any point $Q \in C$ the first
   cohomology group $H^1 \bigl(C[1], F^*(E)\otimes I_B \otimes I_Q
   \bigr)$ will vanish if $b_p(\delta)+1\geq\#(B)+1$, which was one of our
   hypotheses.  The global generation then follows trivially from the
   long exact sequence associated with the following ideal sheaf
   sequence
   $$
   0 \to F^*(E) \otimes I_B \otimes I_Q \to F^*(E)\otimes I_B \to (F^*(E)\otimes I_B)|_Q \to 0.
   $$
   This ends the proof.
 \end{proof}

 We will later need that the vanishing statement of
 Proposition~\ref{prop:vanishing} is an open property for families of
 vector bundles --- recall that a family $(E_t)_{t \in T}$ of vector
 bundles on $C$ is given by an algebraic variety $T$ and a vector
 bundle $E$ on the product $C \times T$. Using that the vanishing
 statement involves all subschemes $B \subset C$ of bounded length,
 this follows from the properness of the Hilbert-scheme by standard
 arguments.

 \begin{lem}[Vanishing after pull-back via Frobenius in families]\label{lem:atopfreefam}
   Let $(E_t)_{t \in T}$ be an irreducible family of vector bundles on
   $C$. Assume that there exists a point $0 \in T$ and a number $n \in
   \mathbb N$, such that for all subschemes $B \subset C$ of length
   $n$, we have $h^1(C, E_0 \otimes I_B) = 0$.

   Then there exists an open neighborhood $T_0 \subset T$ of $0$, such
   that for all $t \in T_0$, and all subschemes $B \subset C$ with
   $\length(B) = n$, we have $h^1(C, E_t \otimes I_B) = 0$. \qed
 \end{lem}

 \section{Relative deformations of curves}
 \label{sec:deformation}
 
 The aim of this section is to fix notation and to gather some
 definitions and results about relative deformations of morphisms from
 curves. While all of this is fairly standard, the relative case is
 not covered in great detail in the standard references. We have
 therefore chosen to say a few words here.

 The setup we are considering is that of a surjective morphism $\sigma:
 C \to X$ of $Y$-schemes,
 $$
 \xymatrix{ C \ar[r]^{\sigma} & X \ar[d]^{\pi} \\ & Y.}
 $$
 Throughout this section we work over an algebraically closed field $k$
 of arbitrary characteristic.

 The main definitions concern the space of $Y$-morphisms.

 \begin{notation}
   Let $\Hom_Y(C, X) \subset \Hom(C, X)$ be the scheme of $Y$-morphisms
   $\sigma : C \to X$. We use $[\sigma]$ to denote the point associated
   with $\sigma$.
 \end{notation}

 \begin{notation}
   Given any such morphism $\sigma$ and a subscheme $B \subset C$, we
   consider the subscheme $\Hom_Y(C,X, \sigma|_B) \subset \Hom_Y(C, X)$
   of $Y$-morphisms $\sigma'$ whose restriction to $B$ agrees with
   $\sigma$, i.e.~$\sigma'|_B = \sigma|_B$. To shorten notation, let
   $\mathcal H_{\sigma, B} \subset \Hom_Y(C, X, \sigma|_B)$ be the
   connected component that contains $[\sigma]$.
 \end{notation}

 When discussing separably uniruled smooth varieties, one of the most
 important notions is that of a ``free morphism''. When the domain is a
 curve, the main feature of a free morphism is that it allows a large
 number of deformations. There is an obvious relative version of this
 notion with similar properties that is not discussed much in the
 literature.

 \begin{defn}
   A $Y$-morphism $\sigma: C \to X$ is called \emph{relatively free}
   over $\sigma|_B$, if
   \begin{enumerate}
   \item $X$ and $\pi$ is smooth along the image of $\sigma$,
   \item $H^1(C, \sigma^*(T_{X|Y}) \otimes I_B) = 0$, and
   \item $\sigma^*(T_{X|Y}) \otimes I_B$ is globally generated.
   \end{enumerate}
 \end{defn}

 \begin{prop}\label{prop:geomfreedom}
   Let $\sigma$ be a $Y$-morphism which is relatively free over
   $\sigma|_B$. If $\dim Y = 1$ and $\pi\circ \sigma$ is surjective,
   then
   \begin{enumerate}
   \item the subscheme $\mathcal H_{\sigma, B}$ is smooth at $\sigma$,

   \item \ilabel{il:geomfreedom2} the universal morphism $\mu :
     C \times \mathcal H_{\sigma, B} \to X$ dominates $X$, and

   \item \ilabel{il:geomfreedom3} if $M\subset X$ is any subset of
     $\codim_X M \geq 2$, and $[\sigma'] \in \mathcal H_{\sigma, B}$ a
     general point, then $(\sigma')^{-1}(M) \subset B$.
   \end{enumerate}
 \end{prop}
 \begin{proof}
   This result is completely similar to \cite[II.~props.~3.5 and
   3.7]{K96} where the non-relative case is considered. The proof is
   analogous.
 \end{proof}

 We close this section by recalling Mori's famous Bend-and-Break
 argument which we use in our proof of Theorem~\ref{thm:algleafs}. This
 result works for absolute and relative deformations alike.

 \begin{thm}[Bend-and-Break, \protect{\cite[prop.~3.3]{K91}}]\label{bendbreak}
   \label{thm:bnb}
   Assume that $C$ is smooth, not rational, that $\sigma$ is not
   constant and that there exists a number $b \in \mathbb N$ such that
   for every closed subscheme $B\subset C$ of length $b$,
   $\dim_{[\sigma]} \Hom(C, X, \sigma|_B) > 0$. If $x \in C$ is a
   general point, and $H$ any ample divisor on $X$, then there exists a
   rational curve $R \subset X$ that contains $\sigma(x)$ and
   satisfies:
   $$
   H\cdot R\leq 2\dfrac{\deg f^*H}{b}.
   $$
 \end{thm}

 \section{Proof of Theorem~\ref*{thm:algleafs}}
 \label{sec:proofOf1}
 
 We maintain the assumptions and notation from
 Theorem~\ref{thm:algleafs} throughout the present section. We will
 work over the complex number field, except in
 Sections~\ref{sec:44}--\ref{sec:45}, where we consider the reduction
 to fields of positive characteristic. Varieties that are defined over
 fields $k$ of positive characteristic will be marked with an index
 $k$.
 
 In Section~\ref{sec:red1} we start by reducing ourselves to the
 situation where the curve $C$ is smooth and everywhere transversal to
 the given foliation $\sF$. We will then, in Section~\ref{sec:algebr}
 recall an old result of Hartshorne which immediately shows that all
 leaves through points of $C$ are algebraic. This argument also reduces
 us to the case where the foliation is the foliation of a morphism that
 is smooth along $C$.
 
 The main point is then to show that the leaves are rationally
 connected. This is done by a reduction to the absurd: we assume in
 Section~\ref{sec:setnot} that there is a non-trivial rationally
 connected quotient. Using Mori's reduction modulo $p$ argument, and
 the vanishing result of Proposition~\ref{prop:vanishing}, we can show
 in Sections~\ref{sec:44}--\ref{sec:45} that the maximally rationally
 connected quotient is covered by rational curves. This contradicts a
 famous result of Graber-Harris-Starr.

 As announced in the introduction, we end with a short section on the
 history of this problem, and with attributions.

 \subsection{Reduction to the case of a normal foliation}
 \label{sec:red1}

\begin{figure}[tbp]
   \centering
     \begin{picture}(11.5,7)(0,0)
      \put(0, 3){\includegraphics[height=4cm]{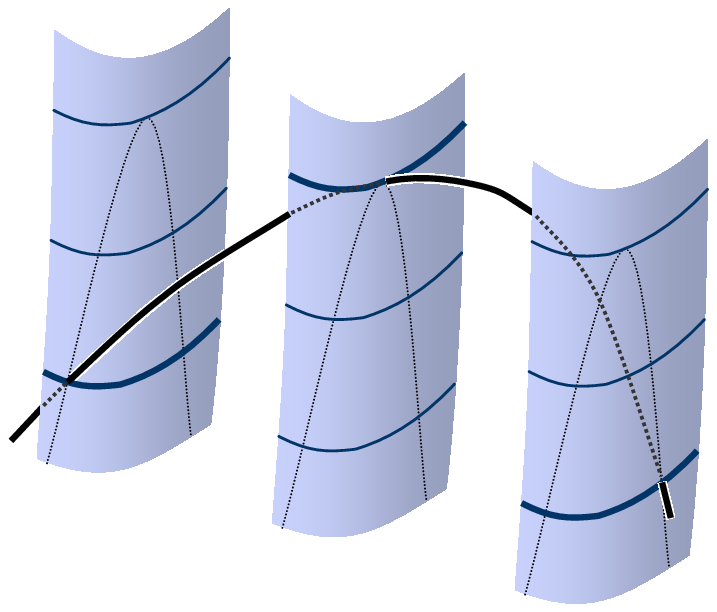}}
      \put(2,1.5){\vector(0,1){1.5}}
      \put(0.7,2.2){${\scriptstyle \sigma :=(\nu,id)}$}
      \put(2.5,3){\vector(0,-1){1.5}}
      \put(2.7,2.2){${\scriptstyle p_2=:\pi}$}
      \put(5.5,5){\vector(1,0){2}}
      \put(6.2,5.2){${\scriptstyle p_1}$}
      \put(8, 3){\includegraphics[height=3cm]{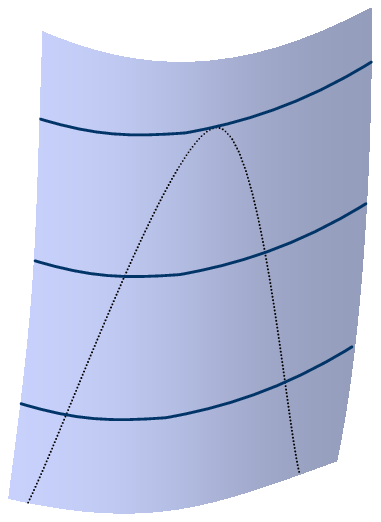}}
      \put(1.5,0.8){$\tilde{C}$}
      \put(2.5,6.7){$Y=X\times\tilde{C}$}
      \put(9,5.9){$X$}
      \put(9.45,4.85){${\scriptstyle C}$}
      \put(1.62,5.5){${\scriptstyle C'}$}
      \put(10.5,4.7){\tiny leaves}
      \put(10.4,4.95){\vector(-1,2){0.3}}
      \put(10.4,4.75){\vector(-1,0){0.3}}
      \put(10.4,4.55){\vector(-1,-2){0.3}}
       \put(0.2,0){\includegraphics[width=4cm]{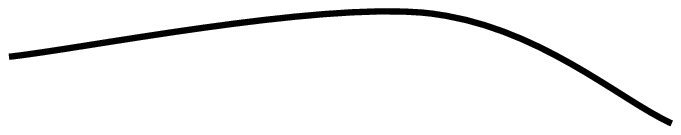}}
     \end{picture}
   \caption{Reduction to the case of a normal foliation}
   \label{fig:normFol}
 \end{figure}

 Let $\nu:\tilde{C}\rightarrow C\subset X$ be a non-constant morphism
 from a smooth curve of positive genus $g(\tilde C) > 0$ to $C$, and
 $Y$ denote the product $X\times \tilde{C}$ with projections $p_1$ and
 $p_2$. We obtain a diagram
 $$
 \xymatrix{
   Y \ar[r]^{p_1} \ar[d]^{p_2=:\pi} & X \\
   \tilde{C} \ar@/^0.3cm/[u]^{\sigma := (\nu,id)} } 
 $$
 Let $C'$ be the image of $\tilde{C}$ under $\sigma$. The sheaf
 $\sF_Y:=p_1^*(\sF)$ is naturally embedded into $T_Y$ via the relative
 tangent bundle $T_{Y|\tilde C}$. It is therefore identified with a
 foliation $\sF_Y$ of $Y$ whose leaves lie in the fibers of $p_2$.
 This construction is depicted in Figure~\ref{fig:normFol}.

 The leaves of $\sF_Y$ are isomorphic to the leaves of $X$, so that it
 suffices to prove the theorem for $\sF_Y\subset T_Y$ and the curve
 $C'$. Clearly the new foliation still fulfills the hypotheses of the
 theorem: it is regular along $C'$ and $\sF_Y|_{C'}\cong \nu^*(\sF|_C)$
 is the pull-back of an ample vector bundle by a finite map, and
 therefore ample.
 
 The construction has the advantage that the curve $C'$ is smooth and
 transversal to $\sF_Y$ since it is already transversal to the
 relative tangent bundle $T_{Y|C}$. Moreover, a leaf of $\sF_Y$
 passing by a point of $C'$ intersects it set--theoretically in
 exactly that point.

 \subsection{The algebraicity of the leaves} 
 \label{sec:algebr}

 Bost and Bogomolov-McQuillan have independently observed that a result
 of Hartshorne on the function fields of formal schemes and the
 ampleness of $\sF_Y|_{C'}$ imply that for any $x \in C'$, the leaf
 through $x$ is algebraic. This answers a question of Miyaoka,
 \cite[rem.~8.9.(1)]{Miy85}. The reader may want to look at
 \cite[thm.~3.5]{Bost01} or \cite[sect.~2.1]{BMcQ} for this. We
 restrict ourselves to a brief review of the argument.
 
 Since $C'$ is everywhere transversal to $\sF_Y$, the classical
 Frobenius theorem, see e.g.~\cite[thm.~1.60]{War71}, immediately
 yields the following.

 \begin{fact}\label{fact:frobstripe}
   There exists an irreducible analytic submanifold $W\subset Y$
   containing $C'$ such that the restriction $\pi|_W$ is smooth, and
   such that its fibers are analytic open subsets of the leaves of the
   foliation passing through points of $C'$. Moreover,
   $\sF_Y|_{C'}\cong N_{C',W} \cong T_{W|\tilde C}|_{C'}$. \qed
 \end{fact}

 Hartshorne's result, applied to a formal neighborhood of $C'$ in $W$,
 then gives the exact dimension of the Zariski closure.

 \begin{fact}[\protect{\cite[thm.~6.7]{Ha68}}]\label{fact:hartshorne}
   Let $\overline{W}$ be the Zariski closure of $W$.  Since $N_{C',W}$
   is ample we have $\dim(\overline{W})=\dim(W)$.  \qed
 \end{fact}

 \begin{cor}\label{cor:algebraicity}
   If $x \in C'$ is any point, and $L \subset (\pi|_{\overline
     W})^{-1}\bigl(\pi(x)\bigr)$ the unique irreducible fiber
   component that contains $W_x := (\pi|_W)^{-1}\bigl(\pi(x)\bigr)$,
   then $L$ is exactly the leaf of $\sF_Y$ through $x$.  In
   particular, all leaves of $\sF_Y$ through points of $C'$ are
   algebraic.
 \end{cor}
 \begin{proof}
   Fact~\ref{fact:hartshorne} immediately implies that $\dim L = \rank
   \sF_Y$. Since $L$ contains the $\sF_Y$-invariant submanifold $W_x$,
   it is clear that $L$ itself is invariant under $\sF_Y$ wherever it
   is smooth.  Hence the claim.
 \end{proof}

 By construction, Corollary~\ref{cor:algebraicity} immediately implies
 that the leaves of $\sF$ through points of $C$ are algebraic. We will
 show in the rest of the present section that they are rationally
 connected, or equivalently, that the fibers of $\overline W \to \tilde
 C$ are rationally connected.

 One technical problem with this approach is that $\overline W$ is not
 necessarily smooth along $C'$. In fact, it may happen that whole
 fibers of $\overline W$ over $\tilde C$ are contained in the
 non-normal locus.  To overcome this difficulty, observe that images of
 rationally connected varieties are themselves rationally connected. It
 suffices therefore to prove rational connectedness for the fibers of
 the normalization of $\overline W$ over $\tilde C$.

 \begin{rem}
   The universal property of the normalization immediately yields an
   embedding $$
   e : W \to \text{normalization of } \overline W.  $$
   In
   particular, the normalization of $\overline W$ is smooth along the
   image $e\circ \sigma (\tilde C)$, and the normal bundle of $e\circ
   \sigma (\tilde C)$ is isomorphic to $\sF_Y|_{C'}$.
 \end{rem}

 \subsection{Rational connectedness of the leaves: setup of notation}
 \label{sec:setnot}
 
 By the previous subsections, we may replace $X$ without loss of
 generality by a desingularization of the normalization of
 $\overline{W}$.  Since we work in characteristic 0, Hartshorne's
 characterization of ampleness \cite[II,~6.4.15]{Laz04} implies that
 if $E$ is any ample vector bundle on $C$, then
 $E\langle\frac{-1}{\rank E +1}\rangle$ is also ample. It therefore
 suffices to prove Theorem~\ref{thm:algleafs} under the following
 additional assumptions that we maintain through the end of the
 present section~\ref{sec:proofOf1}.

 \begin{assumption}\label{ass:fiber}
   In the setup of Theorem~\ref{thm:algleafs}, assume additionally that
   $X$ is smooth, that $C$ is a smooth curve of genus $g(C)>0$, and
   that the foliation $\sF$ is the foliation of a morphism $\pi: X \to
   C$ with connected fibers which is smooth about $C$ and with
   $\pi|_C=\id_C$.

   In particular, if $\sigma : C \to X$ is the inclusion, we assume
   that $\sigma^*(T_{X|C})$ and $\sigma^*(T_{X|C})\langle -\frac{1}{\dim
   X}\rangle$ are ample vector bundles on $C$.
 \end{assumption}
 
 We will need to consider a maximally rationally connected fibration
 of $X$ (``MRC-fibration'') $q : X \dasharrow Z$ which maps $X$ to a
 normal, projective variety $Z$; this is explained in detail in
 \cite[sect.~IV.5]{K96} or \cite[sect.~5]{Debarre01}. The
 functoriality of the MRC-fibration, \cite[thm.~IV.5.5]{K96}, and the
 assumption that $C$ is not rational, imply that $\pi$ factors via the
 morphism $q$. We obtain a diagram as follows,
 \begin{equation}\label{eq:extdiag}
   \xymatrix{
     X \ar[d]^{\pi} \ar@{-->}[r]^{q} & Z \ar@{-->}@/^0.3cm/[ld]^{\beta} \\
      C \ar@/^0.3cm/[u]^{\sigma}}
 \end{equation}
 where all morphisms and maps but $\sigma$ are surjective,
 or dominant, respectively.
 
 We observe that with these notations and assumptions, in order to
 prove Theorem~\ref{thm:algleafs}, it is enough to show:

 \begin{prop}\label{prop:Zcurve}
   The maximally rationally connected quotient $Z$ is a curve.
 \end{prop}
 
 Recall \cite[thm.~5.13]{Debarre01} that the MRC-fibration $q$ is an
 almost holomorphic map. Thus, if Proposition~\ref{prop:Zcurve} holds
 true, then $q$ must be a morphism. Since both $q$ and $\pi$ have
 connected fibers, $\beta$ will be an isomorphism. This will complete
 the proof of Theorem~\ref{thm:algleafs}.
 
 The remainder of section~\ref{sec:proofOf1} is concerned with the
 proof of Proposition~\ref{prop:Zcurve}. For that, we assume the
 contrary, i.e.~that $\dim Z \geq 2$. We will then show that $Z$ is
 uniruled, which contradicts the well-known result of Graber, Harris
 and Starr, \cite{GHS03}. The uniruledness of $Z$ will be established
 by a standard bend-and-break method in positive characteristic that
 is detailed below. The method requires us to bound the degrees of
 curves in $Z$. For this, we have to fix a polarization on $Z$.

 \begin{notation}\label{not:Dandd}
   Choose a point $z \in Z$.  Further, choose a very ample line bundle
   $H_Z$ on $Z$. Let $H_X$ denote the pull-back of $H_Z$ by $q$; the
   rational map $q$ is then determined by a linear system $L\subset
   H^0(X,H_X)$, whose fixed locus is the indeterminacy locus of $q$.
 \end{notation}

 \begin{rem}\label{rem:problemsinsecondasterisque}
   The rational map $\beta$ induces a foliation on $Z$ whose leaves
   are fibers of $\beta$. We will implicitly show that there are lots
   of rational curves in $Z$ that are tangent to this foliation. This
   part of our proof is modelled on \cite[sect.~9]{SecondAsterisque},
   where a similar approach is employed in order to prove
   Corollary~\ref{cor:ratconnleaf}. While in our setup the existence
   of algebraic leaves on $X$ immediately yields the existence of
   Diagram~\ref{eq:extdiag}, the construction of a foliation on the
   partial rationally connected quotient in
   \cite[p.~111]{SecondAsterisque} is problematic: the polarization of
   the $Z$-variety $Q$ discussed in \cite[p.~111]{SecondAsterisque}
   need not be stable under the action of the Galois group of $Q$ over
   $Z$. Accordingly, there is no obvious reason why the sheaf
   $\mathcal L$ is. Also, it is not quite clear in
   \cite{SecondAsterisque} why the images of the $C_t$ are again
   general complete intersection curves.

   At this point of the argumentation, Bogomolov-McQuillan
   \cite[p.~22]{BMcQ} do not consider the MRC-fibration, but a rational
   map whose fibers are single rational curves. We had difficulties
   understanding that part of their paper.
 \end{rem}

 \subsection{Rational connectedness of the leaves: Reduction modulo $p$}
 \label{sec:44}

 In order to prove Proposition~\ref{prop:Zcurve}, we use Mori's
 standard reduction mod $p$ argument, see \cite{Mori79},
 \cite[II.5.10]{K96}. It is then enough to prove
 Proposition~\ref{prop:Zcurve} assuming that all varieties and
 morphisms are defined over an algebraically closed field of large
 characteristic.
 
 To be more precise, let $S$ be a ring of definition of all varieties
 and morphisms that appear in Diagram~\eqref{eq:extdiag}, and of $z$,
 $H_X$ and $H_Z$. Let $X_S$ be the associated scheme over $\Spec S$.
 Given an algebraically closed field $k$ and a morphism $S \to k$, set
 $X_k := X_S \times_{\Spec S} \Spec k$.  We use the same notation for
 $Z_k$, $H_{Z,k}$, etc.
 
 The following Proposition will be shown in Section~\ref{sec:45}
 below.

 \begin{prop}\label{prop:charp}
   The integer
   $$
   d := 2\deg \sigma^*(H_X) \cdot \dim X
   $$
   is positive. There exists a number $p_0$ with the following
   property: If $k$ is the algebraic closure of a residue field of $S$
   such that:
   \begin{itemize}
   \item all assumptions made in subsection~\ref{sec:setnot} still
     hold for the geometric fibers over $k$,
   \item the characteristic of $k$ is larger than $p_0$, i.e. $\charac
     k > p_0$,
   \end{itemize}
   then the geometric fiber $Z_k$ is uniruled with curves of
   $H_{Z,k}$-degrees at most $d$.
 \end{prop}
 
 If Proposition~\ref{prop:charp} holds true, Mori's reduction argument
 implies that $Z = Z_{\mathbb C}$ is uniruled with curves of
 $H_Z$-degree at most $d$. This will complete the proof of
 Proposition~\ref{prop:Zcurve} and hence of
 Theorem~\ref{thm:algleafs}.

 \subsection{Rational connectedness of the leaves: Proof in characteristic $p$}
 \label{sec:45}
 
 We will now prove Proposition~\ref{prop:charp}. Recall from
 Remark~\ref{rem:atft_gets_pos} that there exists a number $p_0$ such
 that for all algebraically closed fields $k$ with $p := \charac(k) >
 p_0$, the vanishing threshold $b_p$, defined on
 page~\pageref{def:atopft} satisfies
 \begin{equation}
   \label{eq:limit2}
   b_p(1/\dim X) > 1 \text{\quad and \quad}
   \left| \, \frac{2p \cdot \deg \sigma^*(H_X)}{\lfloor b_p(1/\dim X)
       \rfloor} - d \, \right| \leq \frac{1}{2}.
 \end{equation}
 
 Now let $k$ be the algebraic closure of a residue field of $S$ that
 satisfies the assumptions of Proposition~\ref{prop:charp}.  As
 before, let $F: C_k[1]\longrightarrow C_k$ be the $k$-linear
 Frobenius morphism.

 \begin{lem}\label{domina}
   There is an open neighborhood $\Omega \subset \Hom_{ C_k}(C_k[1], X_k)$
   of $\sigma_k \circ F$ such that
   \begin{enumerate}
   \item \ilabel{il:domina1} If $[\sigma'] \in \Omega$ is any morphism
     and $B \subset C_k[1]$ any subscheme of length $\#(B) \leq
     b_p(1/\dim X)$, then $\sigma'$ is relatively free over
     $\sigma'|_B$.
   \item \ilabel{il:domina2} If $T \subset X$ is the indeterminacy
     locus of the rational map $q$ then the subset
     $$
     \Omega^0 = \{ [\sigma'] \in \Omega \,\, | \,\,
     (\sigma')^{-1}(T_k) = \emptyset\}
     $$
     of morphisms whose images avoid $T_k$ is again open in
     $\Hom_{C_k}(C_k[1], X_k)$.
   \end{enumerate}
 \end{lem}

 \begin{proof}
   The vanishing results of Proposition~\ref{prop:vanishing} and
   Lemma~\ref{lem:atopfreefam}, applied to $E=\sigma_k^*
   (T_{X_k|C_k})$ and $\delta = \frac{1}{\dim X}$, yield the existence
   of an open set $\Omega$ such that all $[\sigma'] \in \Omega$ are
   relatively free over $\sigma'|_B$ if $\#(B) \leq b_p(1/\dim X)$.
   This shows \iref{il:domina1}.  Assertion \iref{il:domina2} follows
   from Proposition~\ref{prop:geomfreedom}.\iref{il:geomfreedom3}.
 \end{proof}

 \begin{notation}
   Since morphisms $[\sigma'] \in \Omega^0$ avoid $T$, we have a natural morphism:
   $$
   \begin{array}{rccc}
     \eta : & \Omega^0 & \to & \Hom( C_k[1], Z_k ) \\
     & [\sigma'] & \mapsto & [q_k \circ \sigma']
   \end{array}
   $$
   We will also need to consider the associated evaluation morphism.
   $$
   \begin{array}{rccc}
     \mu : & C_k[1] \times \Omega^0 & \to & Z_k \\
     & (y, [\sigma']) & \mapsto & q_k \circ \sigma'(y)
   \end{array}
   $$
 \end{notation}

 The most important properties of $\eta$ and $\mu$ are summarized in
 the following corollary to Lemma~\ref{domina}.

 \begin{cor}\label{cor:morphtoz}
   With the notation introduced above, we have the following:
   \begin{enumerate}
   \item \ilabel{il:morphtoz1} The evaluation morphism $\mu: \tilde
     C_k[1] \times \Omega^0 \to Z_k$ dominates $Z_k$

   \item \ilabel{il:morphtoz3} For all $[\sigma'] \in \Omega^0$ and
     subschemes $B \subset C_k[1]$ of length $\#B \leq b_p(1/\dim
     X)$, we have
     $$
     \dim_{\eta([\sigma'])} \Hom (C_k[1], Z_k, (q_k\circ sigma'
     )|_B ) \geq 1.
     $$

   \item \ilabel{il:morphtoz2} For all morphisms $\tau$ contained in
     the image $\eta(\Omega^0)$, we have
     $$
     \deg \tau^*(H_{Z,k})= p \cdot \deg \sigma^*(H_X).
     $$ 
   \end{enumerate}
 \end{cor}

 \begin{proof}
   Assertions~\iref{il:morphtoz1} and \iref{il:morphtoz3} are immediate
   consequences of
   Proposition~\ref{prop:geomfreedom}.\iref{il:geomfreedom2}.
   
   For Assertion~\iref{il:morphtoz2}, let $[\tau] \in \eta(\Omega^0)$
   be any element. The morphism $\tau$ can then be written as $\tau =
   q_k \circ \sigma'$, where $\sigma': C_k[1] \to X_k$ is a
   deformation of $\sigma\circ F$. This immediately implies
   \begin{align*}
     \deg \tau^*(H_{Z,k}) & = \deg (q_k\circ \sigma')^*(H_{Z,k})= \\
     & = \deg (\sigma')^*(H_{X,k}) = p \cdot \deg \sigma^*(H_X).
   \end{align*}
 \end{proof}

 Let $Z^0$ be the maximal open set contained in the image of the
 evaluation morphism $Z^0 \subset \mu(C_k[1] \times \Omega^0)$.  By
 Corollary~\ref{cor:morphtoz}.\iref{il:morphtoz1} this is not empty.
 Now, if $z \in Z^0$ is any point, then there exists a morphism $\tau
 :C_k[1] \to Z_k$, $[\tau] \in \eta(\Omega^0)$ whose image contains
 $z$.  Bend-and-Break, Theorem~\ref{bendbreak}, then implies that $z$
 is contained in a rational curve $R \subset Z_k$ of degree
 \begin{align*}
   0 < H_{Z,k}\cdot R & \leq \frac{2 \cdot \deg
     \tau^*(H_{Z,k})}{\lfloor b_p(1/\dim X) \rfloor} & &
   \text{Theorem~\ref{thm:bnb},
     Corollary~\ref{cor:morphtoz}.\iref{il:morphtoz3}}\\
   & = \frac{2p \cdot \deg \sigma^*(H_X)}{\lfloor
     b_p(1/\dim X) \rfloor} & & \text{Corollary~\ref{cor:morphtoz}.\iref{il:morphtoz2}}\\
   \Rightarrow 0 < H_{Z,k}\cdot R & \leq d & & \text{l.h.s. is
     integral, Inequality~\eqref{eq:limit2}}.
 \end{align*}
 In particular, the integer $d$ is positive. This shows the first
 statement of Proposition~\ref{prop:charp}.
 
 By \cite[4(c)]{Bourbaki221}, the scheme $\Hom_d(\P^1_k, Z_k)$ of
 non-constant morphisms $f: \P^1_k \to Z_k$ with $\deg f^*(H_{Z,k})
 \leq d$ is quasi-projective. In particular, it contains only finitely
 many irreducible components. The existence of rational curves of
 degree $\leq d$ through every point of $Z^0$ therefore implies that
 there exists one component $\mathcal H \subset \Hom_d(\P^1_k, Z_k)$
 such that the evaluation morphism
 $$
 \mu : \mathcal H \times \P^1_k \to Z_k
 $$
 is dominant. \cite[IV.~prop.~1.4]{K96} then asserts that $Z_k$ is
 uniruled with curves of $H_{Z,k}$-degree at most $d$. This shows
 Proposition~\ref{prop:charp} and therefore ends the proof of
 Theorem~\ref{thm:algleafs}. \qed

 \subsection{Attributions}\label{sec:att1}
 \label{sec:46}

 The reduction to the case of a normal foliation and the proof of the
 algebraicity of the leaves follow \cite{BMcQ} closely. The setup that
 we use to prove the rational connectedness, however, differs somewhat
 from that of \cite{BMcQ}; see
 Remark~\ref{rem:problemsinsecondasterisque}.
 
 The reduction modulo $p$ that Bogomolov-McQuillan and that we employ
 to produce rational curves on the quotient is of course due to Mori
 \cite{Mori79}. The argumentation in characteristic $p$, is different
 from theirs: using the vanishing result of
 Proposition~\ref{prop:vanishing} rather than \cite[lem.~3.2.1]{BMcQ},
 we can give an explicit bound $d$ for the degree of the curves
 constructed. It is obvious in our construction that the number $d$
 does \emph{not} depend on the characteristic; this is perhaps not so
 clear in \cite{BMcQ}. Another advantage of our approach is that we
 construct the curves directly using the standard Bend-and-Break
 argument, and do not have to deal with the partial Frobenius morphism
 associated with a foliation.

\section{Proof of Theorem~\ref*{thm:reeb}}
\label{sec:reeb}

Theorem~\ref{thm:reeb} is an immediate consequence of the following,
stronger statement.

\begin{thm}\label{thm:reeb2}
  Let $X$ be a complex projective manifold, and $\sF \subset T_X$ a
  foliation. Assume that there exists a compact, rationally connected
  leaf $L$ which does not intersect the singular locus of $\sF$. Then
  all leaves are algebraic. The set
  $$
  V := X \setminus \bigcup_{L' \text{non-compact leaf}} \overline{L'}
  $$
  is a Zariski-open neighborhood of $L$, and the restriction
  $\sF|_V$ is the foliation given by a proper submersion $\pi: V \to
  B$. All fibers of $\pi$ are rationally connected.
\end{thm}

\begin{proof}
  
  Recall the standard fact that rationally connected projective
  manifolds are simply connected, \cite[cor.~4.18]{Debarre01}.  The
  holonomy of the foliation along that leaf $L$ is therefore trivial,
  and Reeb's Stability Theorem \cite[thm.~IV.3]{CL85} asserts that
  there exists a fundamental system of analytic neighborhoods of $L$
  in $X$ that are saturated with respect to $\sF$. The triviality of
  the holonomy implies that one of these neighborhoods, say $V^\circ$,
  admits an analytic slice $B^\circ$ and a proper submersion
  $\pi^\circ:V^\circ \to B^\circ$ that induces the foliation
  $\sF|_{V^\circ}$. By \cite[cor.~2.4]{KMM}, all fibers of $\pi^\circ$
  are rationally connected.
  
  Since the normal bundle of $L$ in $X$ is trivial and $L$ is
  rationally connected, we have $h^1(L, N_{L/X}) = 0$.  The Douady
  space $D(X) \cong \Hilb(X)$ parametrizing compact analytic subspaces
  of $X$ is thus smooth at the point representing $L$, and it follows
  immediately that $B^\circ$ embeds as an analytic open subset in a
  component $D$ of $D(X)$. If $U \subset D \times X$ is the universal
  family and $\pi_1, \pi_2$ the canonical projections, we obtain a
  diagram
  $$
  \xymatrix{
    V^\circ \ar[rrr]^{\text{anal. open inclusion}} \ar[d]_{\txt{\scriptsize $\pi^\circ$ \\ \scriptsize ratl.~conn.~fibers}} & & & U \ar[r]^{\pi_2}  \ar[d]^{\pi_1} & X \\
    B^\circ \ar[rrr]_{\text{anal. open inclusion}} & & & D.  }
  $$ 
  
  If $d \in D$ is any point, then the reduced subvariety $X_d :=
  \bigl(\pi_2(\pi_1^{-1}(d))\bigr)_{\red}$ is of pure dimension $\dim
  X_d = \rank \sF$ and is $\sF$-integral at all points $x \in X_d$
  wherever both $X_d$ and $\sF$ are regular. Consequence: if $x \in
  X_d$ is a regular point of $\sF$, then $X_d$ is smooth at $x$ and
  contains the closure of the leaf through $x$ as an irreducible
  component.  It follows that all leaves of $\sF$ are algebraic.
  
  On the other hand, if $x \in V$ is any point, then the associated
  leaf $L'$ is compact. If $y \in \pi_2^{-1}(x)$ is any point of the
  fiber, then $L' = X_{\pi_1(y)}$. By \cite[IV~3.5.2]{K96}, $L'$ is
  rationally chain connected.  Since $L'$ is smooth, it will then be
  rationally connected.  The holonomy argument from above then shows
  that $D$ is smooth at $\pi_1(y)$, and that $\pi_1$ is submersive in
  a neighborhood of $L'$. In particular, the projection $\pi_2$ is
  birational and isomorphic at $y$.
  
  The exceptional locus $E$ of $\pi_2$ does not intersect $V$ and
  therefore does not dominate $D$. Let $B\subset D\setminus \pi_1(E)$
  be the subset of regular points of $D\setminus \pi_1(E)$ over which
  $\pi_1$ is smooth.  We have seen that $\pi_1^{-1}(B)$ is isomorphic
  to $V$ and $\pi:=\pi_1|_{\pi_1^{-1}(B)}$ verifies the required
  properties.
\end{proof}

\section{Proof of Corollary~\ref*{cor:ratconnleaf}}
\label{sec:cor2}

The proof of Corollary~\ref{cor:ratconnleaf} relies on a number of
facts about the Harder-Narasimhan filtration of vector bundles on
curves, which are possibly known to the experts. For lack of an
adequate reference we have included full proofs in
Section~\ref{sec:51} below. We refer to \cite{Seshadri82} for a
detailed account of semistability and of the Harder-Narasimhan
filtration of vector bundles on curves.

\subsection{Vector bundles over complex curves}
\label{sec:51}

To start, we show that any vector bundle on a smooth curve contains a
maximally ample subbundle.

\begin{prop}\label{prop:HN1-new}
  Let $C$ be a smooth complex-projective curve and $E$ a vector bundle
  on $C$, with Harder-Narasimhan filtration
  $$
  0 = E_0 \subset E_1 \subset \ldots \subset E_r = E
  $$
  and $\mu_i:=\mu (E_i/E_{i-1})$ be the slopes of the
  Harder-Narasimhan quotients.  Suppose that $\mu_1>0$ and let $k :=
  \max\{\, i\, |\, \mu_i>0 \}$.  Then $E_i$ is ample for all $1\leq i
  \leq k$ and every ample subsheaf of $E$ is contained in $E_k$.
\end{prop}

\begin{proof}
  Hartshorne's characterization of ampleness
  \cite[thm.~2.4]{Hartshorne71} says that $E_i$ is ample iff all its
  quotients have positive degree. But the minimal slope of such a
  quotient is $\mu_i$ which is positive for all $1\leq i \leq k$.
  
  Let now $F \subset E$ be any ample subsheaf of $E$ and $j :=
  \min\{\, i\, |\, F\subset E_i, \, 1\leq i\leq r \}$.  We need to
  check that $j\leq k$.  By the definition of $j$ and the ampleness of
  $F$, the image of $F$ in $E_j/E_{j-1}$ has positive slope. The
  semi-stability of $E_j/E_{j-1}$ therefore implies $\mu_j>0$ and
  $j\leq k$.
\end{proof}

Proposition~\ref{prop:HN1-new} says that the first few terms in the
Harder-Narasimhan filtration are ample. The following, related
statement will be used later to construct foliations on certain
manifolds.

\begin{prop}\label{prop:ampleness}
  In the setup of Proposition~\ref{prop:HN1-new}, the vector bundles
  $E_j \otimes \bigl(\factor E.E_i.\bigr)^\vee$ are ample for all $0<j
  \leq i<r$.  In particular, if $E_i$ is any ample term in the
  Harder-Narasimhan Filtration of $E$, then $\Hom\bigl(E_i, \factor
  E.E_i.\bigr)$ and $\Hom\bigl(E_i\otimes E_i, \factor E.E_i.\bigr)$
  are both zero.
\end{prop}

\begin{rem}
  If $X$ is a polarized manifold whose tangent bundle contains a
  subsheaf of positive slope, Proposition~\ref{prop:ampleness} shows
  that the first terms in the Harder-Narasimhan filtration of $T_X$
  are special foliations in the sense of Miyaoka,
  \cite[sect.~8]{Miy85}. By \cite[thm.~8.5]{Miy85}, this already
  implies that $X$ is dominated by rational curves that are tangent to
  these foliations.
\end{rem}

\begin{proof}[Proof of Proposition~\ref{prop:ampleness}]
  As a first step, we show that the vector bundle
  $$
  F_{i,j} := \bigl(\factor E_j.E_{j-1}.\bigr)\otimes \bigl(\factor
  E.E_i.\bigr)^\vee
  $$
  is ample.  Assume not. Then, by Hartshorne's ampleness criterion
  \cite[prop.~2.1(ii)]{Hartshorne71}, there exists a quotient $A$ of
  $F_{i,j}$ of degree $\deg_C A \leq 0$.  Equivalently, there exists a
  non-trivial subbundle
  $$
  \alpha : B \to F_{i,j}^\vee = \bigl(\factor
  E_j.E_{j-1}.\bigr)^\vee\otimes\bigl(\factor E.E_i.\bigr)
  $$
  with $\deg_C B \geq 0$. Replacing $B$ by its maximally
  destabilizing subbundle, if necessary, we can assume without loss of
  generality that $B$ is semistable. In particular, $B$ has
  non-negative slope $\mu(B\bigr) \geq 0$.  On the other hand, we have
  that $\bigl(\factor E_j.E_{j-1}.\bigr)$ is semistable.  The slope of
  the image of the induced morphism
  $$
  B\otimes \bigl(\factor E_j.E_{j-1}.\bigr) \to \bigl(\factor E.E_i.\bigr)
  $$
  will thus be larger than $\mu_{\max}\bigl(\factor E.E_i.\bigr)=
  \mu \bigl(\factor E_{i+1}.E_i.\bigr)$.  This shows that $\alpha$
  must be zero, a contradiction which proves the amplitude of
  $F_{i,j}$.
  
  With this preparation we will now prove
  Proposition~\ref{prop:ampleness} inductively.

  \paragraph{\bf Start of induction: $j=1$}
  
  In this case, the above claim and the statement of
  Proposition~\ref{prop:ampleness} agree.

  \paragraph{\bf Inductive Step}
  
  Assume that $1< j \leq i<r$ and that the statement was already shown
  for $j-1$.  Then consider the sequence
  $$
  0 \to \underbrace{E_{j-1}\otimes \bigl(\factor
    E.E_i.)^\vee}_{\text{ample}} \to E_j \otimes \bigl(\factor E.E_i.)^\vee
  \to \underbrace{(\factor E_j.E_{j-1}. )\otimes(\factor
    E.E_i.)^\vee}_{\text{ample}} \to 0
  $$
  But then also the middle term is ample, which shows
  Proposition~\ref{prop:ampleness}.
\end{proof}

\subsection{Proof of Corollary~\ref{cor:ratconnleaf}}
\label{sec:pcorrcl}

We will show that the sheaf $\sF_C$, which is defined only on the
curve $C$ is contained in a foliation $\sF$ which is regular along $C$
and whose restriction to $C$ is likewise ample.
Corollary~\ref{cor:ratconnleaf} then follows immediately from
Theorem~\ref{thm:algleafs}.

An application of Proposition~\ref{prop:HN1-new} to $E := T_X|_C$
yields the existence of a locally free term $E_i \subset T_X|_C$ in
the Harder-Narasimhan filtration of $T_X|_C$ which contains $\sF_C$
and is ample. The choice of $C$ then guarantees that $E_i$ extends to
a saturated subsheaf $\sF \subset T_X$.  The proof is thus finished if
we show that $\sF$ is a foliation, i.e.~closed under the Lie-bracket.
Equivalently, we need to show that the associated O'Neill-tensor
$$
N : \sF \otimes \sF \to \factor T_X.\sF.
$$
vanishes. By Proposition~\ref{prop:ampleness}, the restriction of
the bundle
$$
Hom\left(\sF \otimes \sF, \factor T_X.\sF.\right) \cong (\sF \otimes
\sF)^{\vee} \otimes \factor T_X.\sF.
$$
to $C$ is anti-ample. Ampleness is an open property,
\cite[cor.~9.6.4]{EGA4-28}, so that the restriction of $\sF$ to
deformations $(C_t)_{t\in T}$ of $C$ stays ample for most $t \in T$.
Since the $C_t$ dominate $X$, the claim follows. This ends the proof
of Corollary~\ref{cor:ratconnleaf}. \qed

\subsection{Attributions}\label{sec:att2}
 
The arguments used to derive Corollary~\ref{cor:ratconnleaf} from
Theorem~\ref{thm:algleafs} were certainly known to experts, and are
implicitly contained in the literature, in particular \cite{Miy85} and
\cite{SecondAsterisque}. We would like to thank Thomas Peternell for
explaining the existence of a maximally ample subbundle to us.

\section{A geometric consequence of unstability}
\label{sec:interpret}

Recall that a complex variety $X$ is called $\Q$-Fano if a
sufficiently high multiple of the anticanonical divisor $-K_X$ is
Cartier and ample. The methods introduced above immediately yield that
$\Q$-Fano varieties whose tangent bundles are unstable allow sequences
of rational maps with rationally connected fibers.

\begin{cor}\label{cor:fgnb}
  Let $X$ be a normal complex $\Q$-Fano variety and $L_1, \ldots,
  L_{\dim X-1} \in \Pic(X)$ be ample line bundles. Let 
  $$
  \{0\} = E_{-1} = E_0 \subset E_1 \subset \cdots \subset E_m = T_X
  $$
  be the Harder-Narasimhan filtration of the tangent sheaf with
  respect to $L_1, \ldots, L_{\dim X-1}$ and set
  $$
  k := \max \{ 0 \leq i \leq m \,|\, \mu(E_i/E_{i-1}) > 0\}.
  $$
  Then $k > 0$, and there exists a commutative diagram of dominant
  rational maps
  \begin{equation}
    \label{eq:ratlquot}
    \xymatrix{
      X \ar@{-->}[d]_{q_1} \ar@{=}[r] & X \ar@{-->}[d]_{q_2} \ar@{=}[r] & \cdots \ar@{=}[r] & X \ar@{-->}[d]_{q_k} \\
      Q_1 \ar@{-->}[r] & Q_2 \ar@{-->}[r] & \cdots \ar@{-->}[r] & Q_k, }
  \end{equation}
  with the following property: if $x \in X$ is a general point, and
  $F_i$ the closure of the $q_i$-fiber through $x$, then $F_i$ is
  rationally connected, and its tangent space at $x$ is exactly $E_i$,
  $T_{F_i}|_x = E_i|_x$.
\end{cor}
\begin{proof}
  Let $C \subset X$ be a general complete intersection curve with
  respect to $L_1, \ldots, L_{\dim X-1}$. Since $c_1(T_X) \cdot C >
  0$, Proposition~\ref{prop:HN1-new} implies $k>0$ and that the
  restrictions $E_1|_C, \ldots, E_k|_C$ are ample vector bundles. We
  have further seen in Section~\ref{sec:pcorrcl} that the $(E_i)_{1
    \leq i \leq k}$ give a sequence of foliations with algebraic and
  rationally connected leaves.
  
  To end the construction of Diagram~\eqref{eq:ratlquot}, let $q_i : X
  \dasharrow \Chow(X)$ be the map that sends a point $x$ to the
  $E_i$-leaf through $x$, and let $Q_i := {\rm Image}(q_i)$.
\end{proof}

\begin{rem}
  Corollary~\ref{cor:fgnb} also holds in the more general setup where
  $X$ is a normal variety whose anti-canonical class is represented by
  a Weil divisor with positive rational coefficients.
\end{rem}

It is of course conjectured that the tangent bundle of a Fano manifold
$X$ with $b_2(X)=1$ is stable. We are therefore interested in a
converse to Corollary~\ref{cor:fgnb} and ask the following.

\begin{question}
  Given a $\Q$-Fano variety and a sequence of rational maps with
  rationally connected fibers as in Diagram~\eqref{eq:ratlquot}, when
  does the diagram come from the unstability of $T_X$ with respect to
  a certain polarization? Is Diagram~\eqref{eq:ratlquot} characterized
  by universal properties?
\end{question}

\begin{question}
  To what extent does Diagram~\eqref{eq:ratlquot} depend on the
  polarization chosen?
\end{question}

\begin{question}
  If $X$ is a uniruled manifold or variety, is there a polarization
  such that the MRC-fibration comes from the Harder-Narasimhan
  filtration of $T_X$?
\end{question}


\begin{thebibliography}{GHS03}

\bibitem[BM01]{BMcQ}
Feodor~A. Bogomolov and Michael~L. McQuillan.
\newblock Rational curves on foliated varieties.
\newblock IHES Preprint, February 2001.

\bibitem[Bos01]{Bost01}
Jean-Beno{\^{\i}}t Bost.
\newblock Algebraic leaves of algebraic foliations over number fields.
\newblock {\em Publ. Math. Inst. Hautes \'Etudes Sci.}, 93:161--221, 2001.

\bibitem[CLN85]{CL85}
C{\'e}sar Camacho and Alcides Lins~Neto.
\newblock {\em Geometric theory of foliations}.
\newblock Birkh\"auser Boston Inc., Boston, MA, 1985.
\newblock Translated from the Portuguese by Sue E. Goodman.

\bibitem[Deb01]{Debarre01}
Olivier Debarre.
\newblock {\em Higher-dimensional algebraic geometry}.
\newblock Universitext. Springer-Verlag, New York, 2001.

\bibitem[Fle84]{Flenner84}
Hubert Flenner.
\newblock Restrictions of semistable bundles on projective varieties.
\newblock {\em Comment. Math. Helv.}, 59(4):635--650, 1984.

\bibitem[GHS03]{GHS03}
Tom Graber, Joe Harris, and Jason Starr.
\newblock Families of rationally connected varieties.
\newblock {\em J. Amer. Math. Soc.}, 16(1):57--67 (electronic), 2003.

\bibitem[Gro66]{EGA4-28}
Alexandre Grothendieck.
\newblock \'{E}l\'ements de g\'eom\'etrie alg\'ebrique. {IV}. \'{E}tude locale
  des sch\'emas et des morphismes de sch\'emas. {III}.
\newblock {\em Inst. Hautes \'Etudes Sci. Publ. Math.}, 28:255, 1966.

\bibitem[Gro95]{Bourbaki221}
Alexandre Grothendieck.
\newblock Techniques de construction et th\'eor\`emes d'existence en
  g\'eom\'etrie alg\'ebrique. {IV}. {L}es sch\'emas de {H}ilbert.
\newblock In {\em S\'eminaire Bourbaki, Vol.\ 6}, pages Exp.\ No.\ 221,
  249--276. Soc. Math. France, Paris, 1995.

\bibitem[Har68]{Ha68}
Robin Hartshorne.
\newblock Cohomological dimension of algebraic varieties.
\newblock {\em Ann. of Math.}, 88(2):403--450, 1968.

\bibitem[Har71]{Hartshorne71}
Robin Hartshorne.
\newblock Ample vector bundles on curves.
\newblock {\em Nagoya Math. J.}, 43:73--89, 1971.

\bibitem[Hör05]{Hoe05}
Andreas Höring.
\newblock Uniruled varieties with splitting tangent bundle.
\newblock preprint math.AG/0505327, 2005.

\bibitem[Kol91]{K91}
J{\'a}nos Koll{\'a}r.
\newblock Extremal rays on smooth threefolds.
\newblock {\em Ann. Sci. \'Ecole Norm. Sup. (4)}, 24(3):339--361, 1991.

\bibitem[Kol92]{SecondAsterisque}
J{\'a}nos Koll{\'a}r, editor.
\newblock {\em Flips and abundance for algebraic threefolds}.
\newblock Soci\'et\'e Math\'ematique de France, Paris, 1992.
\newblock Papers from the Second Summer Seminar on Algebraic Geometry held at
  the University of Utah, Salt Lake City, Utah, August 1991, Ast\'erisque No.
  211 (1992).

\bibitem[Kol96]{K96}
J{\'a}nos Koll{\'a}r.
\newblock {\em Rational curves on algebraic varieties}, volume~32 of {\em
  Ergebnisse der Mathematik und ihrer Grenzgebiete. 3. Folge. A Series of
  Modern Surveys in Mathematics}.
\newblock Springer-Verlag, Berlin, 1996.

\bibitem[KMM92]{KMM}
J{\'a}nos Koll{\'a}r, Yoichi Miyaoka, and Shigefumi Mori.
\newblock Rationally connected varieties.
\newblock {\em J. Algebraic Geom.}, 1(3):429--448, 1992.

\bibitem[Lan04]{Langer04}
Adrian Langer.
\newblock Semistable sheaves in positive characteristic.
\newblock {\em Ann. of Math. (2)}, 159(1):251--276, 2004.

\bibitem[Laz04]{Laz04}
Robert Lazarsfeld.
\newblock {\em Positivity in algebraic geometry}.
\newblock Springer-Verlag, Berlin, 2004.

\bibitem[Miy85]{Miy85}
Yoichi Miyaoka.
\newblock Deformation of a morphism along a foliation.
\newblock In S.~Bloch, editor, {\em Algebraic Geometry}, volume~46 of {\em
  Proceedings of Symposia in pure Mathematics}, pages 245--269, Providence,
  Rhode Island, 1985. American Mathematical Society.

\bibitem[Mor79]{Mori79}
Shigefumi Mori.
\newblock Projective manifolds with ample tangent bundles.
\newblock {\em Ann. of Math. (2)}, 110(3):593--606, 1979.

\bibitem[Ses82]{Seshadri82}
C.~S. Seshadri.
\newblock {\em Fibr\'es vectoriels sur les courbes alg\'ebriques}, volume~96 of
  {\em Ast\'erisque}.
\newblock Soci\'et\'e Math\'ematique de France, Paris, 1982.
\newblock Notes written by J.-M. Drezet from a course at the \'Ecole Normale
  Sup\'erieure, June 1980.

\bibitem[War71]{War71}
F.~Warner.
\newblock {\em Foundations of Differentiable Manifolds and Lie Groups}.
\newblock Scott, Foresman and Company, Glenview, Illinois and London, 1971.
\end{thebibliography}
\end{document}